\documentclass[11pt]{amsart}
\usepackage{amssymb, latexsym}
\newtheorem{defin}{}
\newtheorem{saetze}[defin]{}
\newtheorem{conjec}[defin]{}
\newtheorem{lemmas}[defin]{}
\newtheorem{folger}[defin]{}
\newtheorem{bemerk}[defin]{}

\newcommand{\fillbox}{\mbox{$\bullet$}}

\newcommand{\N}{\mathbb N}

\newcommand{\Z}{\mathbb Z}
\newcommand{\Q}{\mathbb Q}

\begin{document}
\title{Polycyclic groups: \\ A new platform for cryptology?}

\author[Bettina Eick]{ Bettina Eick }
\address{Bettina Eick, Institut Computational Mathematics, TU Braunschweig,
Pockelsstr. 14, 38106 Braunschweig, Germany}
\email{beick@tu-bs.de}

\author[Delaram Kahrobaei]{ Delaram Kahrobaei }
\address{ Delaram Kahrobaei, Mathematical Institute, University of St
Andrews, North Haugh, St Andrews, Fife KY16 9SS Scotland, UK}
\email{delaram.kahrobaei@st-andrews.ac.uk}
\urladdr{http://www-groups.mcs.st-and.ac.uk/${\sim}$delaram/}

\thanks{}

\date{2004; submitted}


\begin{abstract}
We propose a new cryptosystem based on polycyclic groups. The
cryptosystem is based on the fact that the word problem can be
solved effectively in polycyclic groups, while the known solutions
to the conjugacy problem are far less efficient.
\end{abstract}

\maketitle

\section{Introduction}
\label{intro} Key exchange problems are of central interest in
cryptology. The basic aim is that two people who can only
communicate via an insecure channel want to find a common secret
key. There are many approaches available which try to solve this
problem. The most classical of these is perhaps the
Diffie-Hellmann key exchange.

Key exchange methods are usually based on one-way functions; that
is, functions which are easy to compute, while their inverses are
difficult to determine. Here 'easy' and 'difficult' can mean that
the complexities or the practicality of the methods are far away
from each other; ideally, the one-way function has a polynomial
complexity and its inverse has an exponential complexity.

Many of the known one-way functions have a common problem: it is
often easy to find a one-way function with a polynomial
complexity, but showing that there is no inverse function with
similar complexity or practicality is usually the difficult part
of the project, since the best inverse function might just not
have been discovered yet. Hence it is of interest to investigate
new one-way functions.

Here we propose a new one-way using similar ideas as in the
Arithmetica key exchange \cite{AAG99}; that is, our one-way
function is based on the word problem and the conjugacy problem in
certain non-commutative groups.

The novelty in our approach is that we propose to use polycyclic
groups as a basis for the protocols: These groups are a natural
generalisation of cyclic groups, but they are much more complex in
their structure than cyclic groups. Hence their algorithmic theory
is more difficult and thus it seems promising to investigate
classes of polycyclic groups as candidates to have a more
substantial platform perhaps more secure.

\section{The Diffie-Hellman key exchange}
\label{diffie}
The Diffie-Hellman key exchange is both the
original public-key idea and an important mechanism in current
use. The practical point is that symmetric ciphers are generally
much faster than asymmetric ones (both hardware and software
reasons), so the public-key cipher is merely used to set up a
private key, a session key, intended to be used only for a single
conversation. This protocol is based on the discrete logarithm,
and it uses the relative difficulty of computing discrete
logarithms. Originally, it was based on the ${\mathbb Z}/p$
discrete-logarithm, but it also generalizes to arbitrary finite
fields, elliptic curves, or to any algebraic structure where
logarithms make sense. We recall the classical idea of the
Diffie-Hellman scheme as follows; see \cite{G00}: First, Alice,
and Bob agree on a large prime $m$, and a primitive root $g$
modulo $m$. These need not be kept secret and can be shared by a
group of users. Alice chooses a large random integer $x$,
privately computes $X = g^x \text{ mod } m$, and sends $X$ to Bob
by the possible insecure channel. Meanwhile, Bob similarly chooses
a large random integer $y$, privately computes $Y = g^y \text{ mod
} m$, and sends $Y$ to Alice across the possible insecure channel.
Then Alice privately computes $k= Y^x \text{ mod } m$, and Bob
symmetrically computes $k' = X^y \text{ mod } m$. Then $k=k'$ mod
$m$, as \[ k' = X^y = {(g^x)}^y = g^{x y} = {(g^y)}^x = k \bmod
m,\] and thus $k=k'$ is the common secret of Alice and Bob. But no
one else on the network can determine $k$, unless they can compute
discrete logarithms.

\section{Key exchange based on non-commutative groups}
\label{crypto} In 1999, the Arithmetica key exchange \cite{AAG99}
was introduced by Anshel, Anshel and Goldfeld. It intends to
achieve the same effect as the Diffie-Hellman key exchange; that
is, it establishes a shared secret when the only communication
possible is across an insecure channel. By contrast with
Diffie-Hellman, it is based on combinatorial group theoretic
properties and it uses non-commutative groups such as braid
groups, making use of the difference of the complexity of the word
and the conjugacy problem in such groups; see \cite{BKL98}.

Since then, the area of 'group-theoretic cryptosystems' has been
very active. For example, Bridson and Howie \cite{BH03} have
considered hyperbolic groups as a basis for a group-theoretic
cryptosystem. Also, Gebhardt \cite{VG03} describes a fast
algorithm to check conjugacy in certain braid groups indicating
that braid groups may not provide a secure basis for
cryptosystems. A class of groups which provides a provably secure
basis for the Arithmetica key exchange seems not to be known so
far.

Below we recall two public key exchange methods: the Arithmetica
key exchange \cite{AAG99} and a version of the Diffie-Hellman key
exchange for non-commutative groups. We use these two key exchange
methods as a basis for the remainder of this paper. We refer to
\cite{AAG99} for their application to braid groups.
\newpage
\subsection{Arithmetica key exchange}
Let $G$ be a finitely generated group with solvable word problem.
Let $S$ and $T$ be two finitely generated subgroups of $G$ and let
$\{ s_1, \cdots, s_n \}$ and $\{ t_1,\cdots, t_m \}$ be generating
sets for $S$ and $T$, respectively. Note that here $x^y$ stands
for $y^{-1} x y$. Now suppose two people, Alice and Bob, want to
agree on a key. The group $G$, its subgroups $S$ and $T$ and their
generators are public information. Then
\begin{itemize}
\item[a)] Alice chooses a secret element $a \in S$ as a word in
the generators $a = s_{i_1}^{a_1} \cdots s_{i_l}^{a_l}$ and
publishes $t_1^a, \ldots, t_m^a$.

\item[b)] Bob chooses a secret element $b \in T$ as a word in the
generators $b = t_{i_1}^{b_1} \cdots t_{i_h}^{b_h}$ and publishes
$s_1^b, \ldots, s_n^b$.
\end{itemize}

Based on this setup, Alice and Bob can use $[a,b]$ as a shared
secret. It is straightforward to observe that both, Alice and Bob,
can compute $[a,b]$ readily, since
\begin{eqnarray*}
[a,b] &=& (b^{-1})^a b = ((t_{i_1}^a)^{b_1} \cdots
(t_{i_h}^a)^{b_h})^{-1} \cdot b \;\;\;\; \mbox{ computable for Bob } \\
      &=& a^{-1} a^b = a^{-1} \cdot (s_{i_1}^b)^{a_1} \cdots (s_{i_l}^b)^{a_l}
          \;\;\;\;\;\;\;\;\; \mbox{ computable for Alice }
\end{eqnarray*}
However, to determine $[a,b]$ based on the published data, we have
to compute $a$ and $b$, respectively. These can be determined by
solving the conjugacy problem and finding an element which
conjugates $t_i$ on $t_i^a$ for $1 \leq i \leq m$ and an element
which conjugates $s_j$ on $s_j^b$ for $1 \leq j \leq n$. Thus the
conjugacy problem can be used to break this cryptosystem.

\subsection{Non-commutative Diffie-Hellman key exchange}
Let $G$ be a non-abelian group with solvable word problem. Let $u
\in G$ and let $S$ and $T$ be two subgroups of $G$ such that
$[S,T] = \{1\}$. Suppose that two people, Alice and Bob, want to
agree on a key. The group $G$, its element $u$ and its subgroups
$S$ and $T$ are public information. Then
\begin{itemize}
\item[a)] Alice chooses a secret element $w \in S$ and publishes
$u^w$.

\item[b)] Bob chooses a secret element $v \in T$ and publishes
$u^v$.
\end{itemize}
If $w$ and $v$ commute, then Alice and Bob can use $u^{wv} =
u^{vw}$ as a shared secret. This is straightforward to determine
for Alice and Bob based on their secret data. The determination of
the shared secret is less easy if only the public data is
available. In this case the conjugacy problem can be used to
determine $w$ and $v$ from $u$ and $u^w$ and $u^v$.

\section{Polycyclic groups for cryptosystems}
\label{polycyc} In this section we consider the use of polycyclic
groups as a basis for the key exchange methods as described above.
Recall that a group is called polycyclic if there exists a
polycyclic series through the group; that is, a subnormal series
of finite length with cyclic factors. There are two different
natural representations for these groups which can be used for
computations: polycyclic presentations and matrix groups over the
integers. We consider these two representations in the following.
\subsection{Polycyclic presentations}
Every polycyclic group has a finite presentation which exhibits
the polycyclic structure of the considered group: a polycyclic
presentation of the form
\[ \langle a_1, \ldots, a_n \mid
    a_j^{a_i} = w_{ij}, a_j^{a_i^{-1}} = v_{ij},
    a_k^{r_k} = u_{kk} \mbox{ for } 1 \leq i < j \leq n \mbox{ and }
    k \in I \rangle\]
where $I \subseteq \{1, \ldots, n\}$ and $r_i \in \N$ if $i \in I$
and the right hand sides $w_{ij}, v_{ij}, u_{jj}$ of the relations
are words in the generators $a_{j+1}, \ldots, a_n$. Using
induction, it is straightforward to show that every element in the
group defined by this presentation can be written in the form
$a_1^{e_1} \cdots a_n^{e_n}$ with $e_i \in \Z$ and $0 \leq e_i <
r_i$ if $i \in I$.

A polycyclic presentation is called consistent if every element in
the group defined by the presentation can be represented uniquely
by a word of the form $a_1^{e_1} \cdots a_n^{e_n}$ with $e_i \in
\Z$ and $0 \leq e_i < r_i$ if $i \in I$. In this case these words
are called normal words. We note that every polycyclic group has a
consistent polycyclic presentation and these presentations are
frequently used as a basis for computations with polycyclic
groups. We refer to \cite{S94} for background and a more detailed
introduction to polycyclic presentations.

\subsubsection{The word problem}
The word problem in a consistent polycyclic presentation can be
solved effectively using the so-called collection algorithm, see
\cite{S94}. This algorithm computes the unique normal word for an
arbitrary word in the generators. The basic idea of the collection
algorithm is that it applies iteratedly the power and conjugate
relations of the given presentation to subwords of a given word
and thus it modifies the given word. The nature of the relations
asserts that an iteration of this process will eventually produce
a normal word. The efficiency of the collection algorithm depends
critically on the sequence of chosen subwords which are processed.
There are various strategies which have been investigated for this
purpose. We refer to \cite{LGS90} and \cite{BH04} for an analysis
of strategies in finite polycyclic groups and in \cite{VG03} for
arbitrary polycyclic groups. The resulting complexities of the
methods depend on the growth of the exponents $e_j$ of generators
$g_{i_j}$ occurring in intermediate stages of the algorithm while
processing the word $g_{i_1}^{e_1} \cdots g_{i_l}^{e_l}$. This
growth can be bounded above if the considered group is finite. In
infinite polycyclic groups there is the potential risk of an
integer explosion inherent in the collection algorithm. In praxis,
collection is known as an effective method to solve the word
problem in consistent polycyclic presentations. The method is
implemented in GAP \cite{gap} and MAGMA \cite{magma} and it has
proved to be practical for finite and infinite polycyclic groups.

\subsubsection{The multiple conjugacy problem}
The Arithmetica key exchange can be broken if an effective
algorithm to determine a multiple conjugating element $a$ with
$r_i^a = s_i$ for $1 \leq i \leq m$ can be found (knowing that
such an element exists). We observe that this problem reduces to
the single conjugacy problem as follows.

Let $a_1$ be an element with $r_1^{a_1} = s_1$ and let $G_1 =
C_G(s_1)$. Then every element $a \in G$ with $r_1^a = s_1$ can be
written as $a = a_1 c$ with $c \in G_1$.

By induction, suppose that an element $a_j$ is given with
$r_k^{a_j} = s_k$ for $1 \leq k \leq j$ and let $G_j =
C_{G_{j-1}}(s_j) = C_G(s_1, \ldots, s_j)$. Then every element $a
\in G$ with $r_i^a = s_i$ for all $i$ is of the form $a = a_j c$
with $c \in G_j$. Compute $c \in G_j$ with $(r_{j+1}^{a_j})^c =
s_{j+1}$ and, simultaneously, determine $G_{j+1} =
C_{G_j}(s_{j+1})$. Now define $a_{j+1} = a_j c$ as the element for
the next step. Iterating this process yields an element $a = a_m$
eventually.

As every centralizer $G_j$ is a subgroup of $G$, an induced
polycyclic presentation for $G_j$ can be computed from a
generating set. Thus the multiple conjugacy problem reduces to $m$
applications of the single conjugacy problem with a simultaneous
determination of the corresponding centralizers.

\subsubsection{The single conjugacy problem}

The non-commutative Diffie-Hellmann key exchange and the
Arithmetica key exchange can be broken if an effective algorithm
to determine a conjugating element $r^a = s$ for given $r,s \in G$
can be found. For the Arithmetica key exchange the additional
computation of $C_G(s)$ is necessary.

In \cite{BCRS} has been proved that there exists an algorithm to
compute such a conjugating element and such a centralizer in a
polycyclic group. Later, in \cite{EO03} another algorithm for this
purpose has been described. The algorithm of \cite{EO03} has been
implemented in the Polycyclic package \cite{poly} based on GAP and
Kant \cite{kant} and its performance has been investigated. It can
be observed that the algorithm is practical on interesting
examples, but its performance is far less good than the
performance of the collection algorithm. In particular, in large
polycyclic groups with a complex structure, the computation of
conjugating elements is practically impossible while the word
problem is still effectively solvable.

The algorithm described in \cite{EO03} uses induction down along a
normal series with elementary or free abelian factors of the
considered group $G$. The steps corresponding to finite factors in
the series can be solved by an application of an orbit-stabilizer
algorithm as typical for finite groups.

The steps corresponding to infinite factors require:
\begin{itemize}
\item[a)] Methods from representation theory such as the
computation of submodule series,

\item[b)] Methods from algebraic number theory such as the
computation of unit groups,

\item[c)] Methods from finite polycyclic groups such as the finite
orbit-stabilizer algorithm.
\end{itemize}

A determination of the complexity of the conjugacy algorithm has
not been attempted yet. However, the complexity of this algorithm
could incorporate the complexity of the unit group computation in
algebraic number fields and in this case it is going to be far
away from the complexity of the word problem. This seems to
suggest that polycyclic groups are a promising candidate for
cryptosystems.

\subsection{Matrix representations}
Every polycyclic group can be described as a finitely generated
subgroup of $GL(d, \Z)$ for some $d \in \N$. In this setting, the
word problem for a polycyclic group is solvable in polynomial
time, as matrix multiplication for such groups is solvable in
polynomial time.

The conjugacy problem has not been considered for polycyclic
matrix groups so far. The only implemented method to solve the
multiple or single conjugacy problem in such groups is by
determining a polycyclic presentation and then applying the method
of \cite{EO03}. However, there is a possibility that linear
methods can be used to compute conjugating elements or to improve
the above approach, even though such methods have not been studied
yet.

\section{Examples}
In this section we consider a few examples of polycyclic groups
and we investigate the practicality of the collection algorithm
and the conjugacy algorithm of \cite{EO03}. This will yield the
observation that not every class of polycyclic will be useful as a
basis for a cryptosystem as outlined above, but there are classes
which look promising.

\subsection{Some example computations}
Let $K = \Q[x]/(f)$ be an algebraic number field for a cyclotomic
polynomial $f_w$, where $w$ is a primitive $r$-th root of unity.
Then $deg(f_w) = \varphi(n)$. The maximal order $O$ of $K$ is a
ring whose additive group is isomorphic to $\Z^n$ and the unit
group $U$ of $K$ is a finitely generated abelian group. The
natural split extension $G(w) = O \rtimes U$ is a metabelian
polycyclic group and we investigate a few examples of such groups
and the practicality of their word and conjugacy problem in the
following tables. For this purpose, we list the order $r$ of $w$,
the Hirsch length $h(G(w))$ of $G(w)$, the average time used for
100 applications of the collection algorithm on random words and
the average time used for 100 applications of the conjugacy
algorithm on random conjugates in such groups.

\begin{table}[htb]
\begin{center}
\begin{tabular}{cccc}
r & h(G(w)) &  coll & conj \\
\hline
3 & 2   & 0.00 sec &  9.96 sec \\
4 & 2   & 0.00 sec &  9.37 sec \\
7 & 6   & 0.01 sec & 10.16 sec \\
11 & 14 & 0.05 sec & $>$ 100 hrs
\end{tabular}
\end{center}
\end{table}
For the prime 11, the result of a single conjugacy test could not
be computed within one hour using the current methods. For primes
larger than 11, the results of such timing are expected to be even
worse.

\subsection{Nilpotent groups}
For finitely generated nilpotent groups there are various special
methods known to compute with these groups. For example, the word
problem in finitely generated nilpotent groups can be solved by
evaluating polynomials as described in \cite{LGS98}. The conjugacy
problem can be solved by an efficient methods as for example
described in \cite{S94}.

In this case the complexity and the practicality of the conjugacy
problem is perhaps not as far away from the complexity of the word
problem as in the case of an arbitrary polycyclic group.

\bibliographystyle{abbrv}
\bibliography{XBib}
\end{document}